\documentclass[a4]{amsart}
\usepackage{setspace}
\usepackage{a4}
\usepackage{amssymb,amsmath,amsthm,latexsym}

\usepackage[left=20mm,top=0.6in,bottom=12mm]{geometry}

\usepackage{amsfonts}
\usepackage{amsfonts}
\usepackage{graphicx}
\usepackage{textcomp}
\usepackage{cite}
\usepackage{enumerate}
\usepackage[mathscr]{euscript}
\usepackage{mathtools}
\newtheorem{theorem}{Theorem}[section]

\newtheorem{definition}[theorem]{Definition}
\newtheorem{example}[theorem]{Example}

\newtheorem{proposition}[theorem]{Proposition}
\newtheorem{remark}[theorem]{Remark}
\title{This is the title}
\usepackage{amssymb}
\usepackage{amssymb}
\usepackage{amssymb}
\usepackage{amssymb}
\usepackage{amsmath}
\usepackage{tikz}
\usepackage{hyperref}
\usepackage{enumerate}
\usepackage{mathtools}
\usepackage{amsmath}
\usepackage{tikz}
\usepackage{amssymb}
\usepackage{amsmath}
\usepackage{tikz}
\usepackage{hyperref}
\def\H{\mathcal H}

\def\calC{\mathcal C}
\def\D{\mathcal D}

\begin{document}
\title
{{Closed $ EP $ and Hypo-$ EP $ Operators on Hilbert Spaces}}
\author{P. Sam Johnson}
\address{Department of Mathematical and Computational Sciences, National Institute of Technology Karnataka (NITK), Surathkal, Mangaluru 575 025, India}
\email{sam@nitk.edu.in}

\maketitle

\begin{center}
{\sl Dedicated to Professor C. Ganesa Moorthy on his 60th birthday}
\end{center}

\begin{abstract} 
A bounded linear operator $ A$ on a Hilbert space $ \mathcal  H $ is said to be an $ EP $ (hypo-$ EP $)  operator if ranges of $ A $ and $ A^* $ are equal (range of $ A $ is contained in range of $ A^* $) and $ A $ has a closed range. In this paper, we define $EP$ and hypo-$EP$ operators for densely defined closed linear operators on Hilbert spaces and extend results from bounded operator settings to (possibly unbounded)  closed operator settings.
\end{abstract}



\section{Introduction}

A square matrix $A$ over the complex field $\mathbb{C}$ is said to be an $EP$ matrix  ($ EP $ stands for {\it Equal Projections})  if ranges of $A$ and $A^*$ are equal. Nevertheless the  $EP$ matrix was defined by Schwerdtfeger \cite{Hans1950} in 1950,  it could not get any greater attention until Pearl \cite{Pearl1966(2)} 
gave an interesting characterization of $EP$ matrix through Moore-Penrose inverse: $A$ is an $EP$ matrix if and only if $A^{\dag}A=AA^{\dag}$. 
Campbell and Meyer \cite{Campbell1975} extended the notion of $EP$ matrix to a bounded linear  operator with a closed range defined on a Hilbert space, using the Pearl's characterization.  Itoh \cite{Itoh2005} introduced hypo-$EP$ operator by weakening the Pearl's characterization as $A^\dag A - A A^\dag$ is a positive operator. 

In what follows we will use the term operator to mean a linear operator with domain and range in  (real or complex)  Hilbert spaces.

 $EP$ matrices, bounded $EP$ operators  and bounded hypo-$EP$ operators have been studied by many authors \cite{Hartwig1997, Koliha2007, Pearl1966,Brock1990, Campbell1975, Sam2018,Sam2017}.  In this paper, we extend characterizations and results of bounded $EP$ and hypo-$EP$ operators to (possibly unbounded) closed operators on Hilbert spaces. 
 We recall that an operator on a Hilbert space $\mathcal  H $ is a bounded operator, that is, it maps bounded sets onto bounded sets, if and only if it is continuous.  It follows that unbounded operators are discontinuous (everywhere).  Although existence of unbounded operators on the whole of the Hilbert space $\H$ can be shown (cf. \cite{nair-fa}), unbounded operators of theoretical and practical interests are {\it closed operators} defined on a subspace of $\H$.
Consequently, the specification of the subspace $\mathcal  D $ on which $A$ is defined, called the \textit{domain} of $A$, denoted by $\mathcal D(A)$ is  to be specified, when $A$ is an unbounded operator.  The null space and range of $A$ are denoted by $\mathcal  N(A)$ and $\mathcal R(A)$ respectively.  $M^\perp$ denotes the orthogonal complement of a set $M$ whereas $M\oplus N$ denotes the orthogonal direct sum of subspaces $M$ and $N$.  For the sake of completeness of exposition, we  first begin with the definition of a closed operator.

  \begin{definition}
  	Let $ A$ be an  operator from a Hilbert space $ \mathcal  H $ with domain $ \mathcal  D(A) $ to a Hilbert space $ \mathcal  K $. If the graph of $ A$ defined by 
  	$$ \mathcal  G(A)=\left\{(x,Ax): x\in \mathcal  D(A)\right\} $$ is closed in $ \mathcal  H\times \mathcal  K $, then $ A $ is called a \textbf{closed  operator}. Equivalently, $ A $ is a closed operator if $ x_n\in \mathcal  D(A) $ such that $ x_n\rightarrow x $ and $ Ax_n\rightarrow y$ for some $ x\in\mathcal  H,y\in \mathcal  K $,  then $ x\in \mathcal  D(A) $ and $ Ax=y $. 
  \end{definition}

By the closed graph theorem,  a closed operator with domain $\D(A)$  is a bounded operator if and only if  $\D(A)$ is closed  (cf.  \cite{nair-fa}, Theorem 7.1.2, page 300). In other words, a closed operator defined on all of $\mathcal  H $ is necessarily bounded.  
Given an operator  $ A$ from  $\mathcal  H$ to $\mathcal  K$ with dense domain $ \mathcal  D(A)$, there exists a unique operator $ A^* $ such that 
   \begin{equation}\label{adjoint}
\langle Ax,y\rangle=\langle x,A^*y \rangle \quad \text{for  } x\in \mathcal  D(A),  y\in \mathcal  D(A^*),
   \end{equation}
where 
$$ \mathcal  D(A^*)=\Big\{y\in \mathcal  K: x\mapsto \langle Ax,y\rangle \mbox{ for all } x\in \mathcal  D(A) \mbox{ is continuous}\Big\}.$$
This operator $A^*$ is known as the \textbf{adjoint} of $ A $. The adjoint of any densely defined (not necessarily closed) operator is always closed.    An operator having a dense domain is called a \textbf{densely defined operator}.
  The set of all densely defined closed operators from $ \mathcal  H $ to $ \mathcal  K $ is denoted by $ \mathcal  C(\mathcal  H,\mathcal  K) $ and we write $\mathcal  C(\mathcal  H) = \mathcal  C(\mathcal  H,\mathcal  H)$. We denote the set of bounded operators from $\mathcal  H$ into $\mathcal  K$ by $\mathcal  B(\mathcal  H, \mathcal  K)$ and $\mathcal  B(\mathcal  H) = \mathcal  B(\mathcal  H, \mathcal  H).$  
  We call $ \mathcal  D(A)\cap \mathcal  N(A)^\perp $, the \textbf{carrier} of $ A $ and it is denoted by $ C(A)$.  We note that,  for any $A\in \mathcal  C(\mathcal  H),$ the closure of $C(A)$, that is,  $ \overline{C(A)}$ is $\mathcal  N(A)^\perp$.
  If $A$ and $B$ are in $\mathcal  C(\mathcal  H)$ with $\mathcal  D(A)\subseteq \mathcal  D(B)$ and $Ax=Bx$ for all $x\in \mathcal  D(A)$, then $B$ is called an \textbf{extension} of $A$ and is denoted by $A\subset B$.  Also, we write $A=B$ if $A\subset B$ and $B\subset A$.   Among the unbounded operators, densely defined closed operators are certainly the nicest ones.  For instance, the adjoint $A^*$ of a densely defined closed operator $A$ exists and it is also densely defined with $A=A^{**}$.  However, there are densely defined unbounded operators $A$ with  $\mathcal D(A^*)$ is zero.  
  
In the sequel, we will need the following definitions and results.
  
  \begin{definition}\cite{Rudin1991, sch-book}
  	Let 	$ A\in \mathcal  C(\mathcal  H)$. The operator $ A $ is said to be 
  	\begin{enumerate}
  		\item \textbf{normal} if  $ AA^*=A^*A $, where either side is defined on its standard domain.
  		\item \textbf{symmetric} if  $A\subset A^*$.
  		\item \textbf{self-adjoint} if $ A=A^* $.
  		\item \textbf{positive} if $A$ is self-adjoint and  $ \left<Ax,x\right>\geq 0 $ for all $ x\in \mathcal  D(A) $,  and denote this fact by $A \geq 0$. 
  		  		  	\end{enumerate}
  	  		  	The number $ \gamma (A)=\inf \{\|Ax\|: x\in C(A), \|x\|=1\} $ is called the \textbf{reduced minimum modulus} of $ A $. Moreover, $ \gamma (A)=\gamma(A^*)$.
The operator $|A|:=(A^*A)^{1/2}$ is called the \textbf{modulus of $A$.}  It can be verified that $\mathcal  D(|A|)=\mathcal  D(A)$ and $\mathcal  N(|A|)=\mathcal  N(A)$ and $\overline{\mathcal  R(|A|)}=\overline{\mathcal  R(A^*)}.$
 \end{definition}
 In the above, the square-root of a positive self-adjoint operator $A$ may be  defined using spectral theorem (cf. \cite{yosida}).  
 A simple proof for the existence of a square-root of a positive
unbounded operator was given by Bernau \cite{MR0226433}.

In the following theorem, we state some standard results from operator theory.

  \begin{theorem}\cite{Ben-Israel2003}\label{Ben_thm1.0}
  	Let $ A\in \mathcal  C(\mathcal  H) $. Then the following are true. 
  	\begin{enumerate}
  		\itemsep0em
  		\item $ \mathcal  N(A)=\mathcal  R(A^*)^\perp,\quad \mathcal  N(A^*A)=\mathcal  N(A). $
  		\item $ \mathcal  N(A^*)=\mathcal  R(A)^\perp,\quad \mathcal  N(AA^*)=\mathcal  N(A^*). $
  		\item $ \overline{\mathcal  R(A)} =\mathcal  N(A^*)^\perp,\quad \overline{\mathcal  R(A)}=\overline{\mathcal  R(AA^*)}.$
  		\item $ \overline{\mathcal  R(A^*)} =\mathcal  N(A)^\perp,\quad \overline{\mathcal  R(A^*)}=\overline{\mathcal  R(A^*A)}.$	
  	\end{enumerate}  	
  \end{theorem}

It is well known that, the Moore-Penrose inverse  $A^\dagger$ of a bounded operator on  $\H$  is a closed densely defined operator with 
$\D(A^\dagger) : = \mathcal R(A)\oplus \mathcal  R(A)^\perp$ (cf. \cite{nair-fa}, Theorem 14.3.7). Analogously, the Moore-Penrose inverse $ A^\dag $ of a closed densely defined operator $A$ can be  defined  with 
$ \mathcal  D(A^\dag):=\mathcal  R(A)\oplus\mathcal  R(A)^\perp $  and taking values in $ C(A)$  by associating each $y\in \D(A^\dagger)$ to the unique $A^\dagger y$ such that 
$$AA^\dag y=Qy,$$
  where $ Q  $ is the orthogonal projection of $ \mathcal  H $ onto $ \overline{\mathcal  R(A)}$. It can be seen that  $ \mathcal  N(A^\dag)=\mathcal  R(A)^\perp $ and 
  $$A^\dag Ax=Px \mbox{ for } x\in \mathcal  D(A),$$
where $ P $ is the orthogonal projection of $ \mathcal  H $ onto $\overline{C(A)}$.
Again,  $ A^\dag $  is a closed densely defined operator (cf. \cite{nair-linop, Nashed1976}).

 \begin{theorem}\cite{Nashed1976}\label{thm_unpen}
  	Let $ A\in \mathcal  C(\mathcal  H)$.  The   conditions listed in each of (1) and (2) below characterize  the Moore-Penrose inverse of $ A $.
  	\begin{enumerate}
  		\item\begin{enumerate}
  			\item $ A^\dag AA^\dag y=A^\dag y $ for all $ y\in \mathcal  D(A^\dag) $,
  			\item $ A^\dag Ax=P_{\overline{\mathcal  R(A^\dag)}} x$ for all $ x\in \mathcal  D(A) $,
  			\item  $ AA^\dag y=P_{\overline{\mathcal  R(A)}} y$ for all $ y\in \mathcal  D(A^\dag) $.
  		\end{enumerate}
  		\item \begin{enumerate}
  			\item  $ AA^\dag Ax=Ax $ for all $ x\in \mathcal  D(A) $,
  			\item $ A^\dag AA^\dag y=A^\dag y $ for all $ y\in \mathcal  D(A^\dag) $,
  			\item $ A^\dag A $ and $ AA^\dag $ are symmetric operators.
  		\end{enumerate}
  	  	\end{enumerate}
  \end{theorem}

\begin{theorem}\cite{Ben-Israel2003}\label{prop-dag}
	Let $ A\in \mathcal C(\mathcal  H) $. Then the following assertions hold good.
	\begin{enumerate}
		\item $ A^{\dagger\dagger}=A $.
		\item $ A^{*\dag}=A^{\dagger*} $.
		\item $ \mathcal N(A^{*\dag})=\mathcal N(A) $.
		\item $ A^*A $ and $ A^\dag A^{*\dag} $ are positive and $ (A^*A)^\dagger=  A^\dag A^{*\dag}$.
		\item $ AA^* $ and $ A^{*\dag} A^\dag$  are positive and $ (AA^*)^\dagger=A^{*\dag}A^\dag $.
		\item $  A$ is bounded if and only if $ \mathcal  R(A^\dag) $ is closed.\label{bdd-iff-dag-closed}
	\end{enumerate}
\end{theorem}

In the following result, we give some characterizations of densely defined closed operators to have closed ranges.

  \begin{theorem}\cite{Ben-Israel2003, tarcsay-2012} \label{Ben_thm_1}
  	Let $ A\in \mathcal  C(\mathcal  H) $. Then
  	the following statements are equivalent:
  	\begin{enumerate}
  		\item $ \mathcal  R(A) $ is closed ;
  		\item $ \mathcal  R(A^*) $ is closed ;
  		\item $ \mathcal  R(A^*A) $ is closed ;
  		\item $ \mathcal  R(AA^*) $ is closed ;
  		\item $ \mathcal R(A) =  \mathcal R(AA^*)$;
  		\item $ \mathcal R(A^*) = \mathcal R(A^*A)$;
  		\item $ A|_{C(A)}$ has a bounded inverse from its range;
  		\item $ A^\dag $ is a bounded operator with $\mathcal D(A^\dag)=\mathcal H$ ;
  		\item $\gamma (A)>0$ ;
  		\item there exists $k>0$ such that $\|Ax\|\geq k \|x\|$ for all $x\in C(A)$ ;   		
  		\item there exists $k>0$ such that $\|A^*x\|\leq k \|AA^*x\|$ for all $x\in \mathcal D(AA^*)$\label{char-11} ;   
  		\item there exists $S\in \mathcal B(\mathcal H)$ such that $A\subset AA^*S.$\label{char-12}
  	\end{enumerate}
  \end{theorem}
  
In the following we list some more properties of operators in $\calC(\H)$.
  
\begin{theorem}\cite{tarcsay-2012}\label{tarcsay}
	Let $ A\in \mathcal  C(\mathcal  H) $. Then the following assertions hold good.
	 \begin{enumerate}
		\item If $\mathcal D(A^*A)=\mathcal D(AA^*)$, then $\mathcal D(A)=\mathcal D(A^*)$.
		\item If $\mathcal D(A)\subset \mathcal D(A^*)$, then $A$ is bounded if and only if $\mathcal D(A)=\mathcal D(A^2)$.
		\item $A$ is bounded if and only if $\mathcal D(AA^*)=\mathcal D(A^*)$.
		\item $A^*$ is bounded if and only if $\mathcal D(A^*A)=\mathcal D(A)$.
	\end{enumerate}
\end{theorem} 

In the following result,  we list some relationships among the concepts, majorization,
  range inclusion, and factorization, which  are studied in a general setting for densely defined closed operators by Douglas in 1966 \cite{Douglas1966}.  Note that (\ref{char-11}) and (\ref{char-12}) in Theorem \ref{Ben_thm_1} give majorization and factorization type characterizations of closed range operators. 
  
 \begin{theorem}\cite{Douglas1966}
 	Let $ A $ and $ B $ be   in  $ \mathcal C(\mathcal H) $. Then the following are true.
 	\begin{enumerate}
 		\item If $ AA^*\leq BB^*$, then there exists a contraction $ C $ so that $ A\subset BC $. 
 		\item If $ C $ is an operator so that $ A\subset BC $, then $ \mathcal R(A) \subseteq \mathcal R(B) $.
 		\item If $ \mathcal R(A) \subseteq \mathcal R(B) $, then there exists a densely defined operator $ C $ so that $ A = BC $ and a
  		number $ k > 0 $ so that
  		$ \|Cx\|^2\leq k\left\{\|x\|^2+\|Ax\|^2\right\} $ for $ x\in \mathcal D(C) $.  Moreover, if $A$ is bounded, then $C$ is bounded ; if $B$ is bounded, then $C$ is closed.
  	\end{enumerate}

  \end{theorem}

Bounded $ EP $ and hypo-$ EP $ operators have been studied by many authors. However,  $ EP $ and hypo-$ EP $ operators,   in the setting of  unbounded operators, have not got much attention  in the literature.  Motivated by the results in \cite{Campbell1975, Brock1990, Koliha2007} for bounded $EP$ operators and in  \cite{Sam2017, Sam2018} for bounded hypo-$EP$ operators, we have made an effort in the paper to glean the results for the class of (possibly unbounded) closed densely defined operators in  Hilbert spaces.

 \section{Closed $ EP $ Operators on Hilbert Spaces}

 We begin the section with the definition of bounded $EP$ operator and some of its characterizations. 
 The notion of $EP$ operator was introduced by Campbell and  Meyer \cite{Campbell1975} in 1975. Brock \cite{Brock1990} gave few more characterizations of $EP$ operators.
 \begin{definition}\cite{Campbell1975}
 	An operator $A\in \mathcal  B(\mathcal  H)$  is called an \textbf{$EP$ operator} 
  if $ A $ has a closed range and  $\mathcal  R(A)=\mathcal  R(A^*)$.
 \end{definition}

As already mentioned in the Introduction, $ EP $ stands for {\it Equal Projections}. 
Note that the set of all bounded $EP$ operators contains all bounded normal operators with closed ranges.
 
 \begin{theorem}\cite{Brock1990}\label{bdd-EP}
 	Let $A\in \mathcal B(\mathcal  H)$ with a closed range. Then the following are equivalent:
 	\begin{enumerate}[(i)]
 		\itemsep0em
 		\item $A$ is an $EP$ operator ;
 		\item $AA^\dag=A^\dag A$ ;
 		\item $\mathcal  N(A)^\bot=\mathcal  R(A)$ ;
 		\item $\mathcal  N(A)=\mathcal  N(A^*)$ ;
 		\item $A^*=PA$, where $P$ is a bijective bounded operator on $\mathcal  H$.
 	\end{enumerate}
 \end{theorem}

The following result characterizes densely defined closed operators with closed ranges to have same ranges of `operators' and `their adjoints'.

 \begin{theorem}\label{ch_6_th_4}
 	Let $ A\in \mathcal  C(\mathcal  H) $ with a closed range. Then the following are equivalent:
 	\begin{enumerate}
 		\item $\mathcal R(A)=\mathcal R(A^*)$ ;\label{un_1}
 		\item $AA^\dag=A^\dag A$ on $ \mathcal  D(A) $ ;\label{un_2}
 		\item $\mathcal  N(A)=\mathcal  N(A^\dag)$ ;\label{un_3}
 		\item $\mathcal  N(A)=\mathcal  N(A^*)$ ;\label{un_4}
 		\item $\mathcal  N(A)^\perp =\mathcal  R(A)$ ; \label{un_7}
 		\item  $ \overline{C(A)}=\mathcal  R(A) $ ;\label{un_5}
 		\item $ \mathcal  H= \mathcal  R(A)\oplus \mathcal  N(A)$.\label{un_6}

 	\end{enumerate}
 \end{theorem}
 \begin{proof}
Suppose that $\mathcal R(A)=\mathcal R(A^*)$. Let $ x\in \mathcal  D(A)=\mathcal  N(A)\oplus C(A) $. Then 
$ x=x_1+x_2, x_1\in \mathcal  N(A), x_2\in C(A)=\mathcal  R(A^\dag) $. Hence $ A^\dag Ax= A^\dag A(x_1+x_2)=A^\dag Ax_2=x_2$. As $ C(A)=\mathcal  R(A^\dag) \subseteq \mathcal  R(A^*)=\mathcal  R(A)$ and $ \mathcal  N(A) =\mathcal  R(A^*)^\perp =\mathcal  R(A)^\perp =\mathcal  N(A^\dag)$.  Now $ AA^\dag x=AA^\dagger(x_1+x_2)=AA^\dagger x_1+AA^\dagger x_2=AA^\dagger x_2=x_2 $. Therefore $ AA^\dag= A^\dag A $ on $ \mathcal  D(A) $.  	Now assume  $ A^\dag A=AA^\dag  $ on $ \mathcal  D(A) $. Then $ \overline{\mathcal  R(A^\dag)}=\overline{\mathcal  R(A)} $ and hence $ \mathcal  R(A^*)=\mathcal  R(A) $. Thus the implication (\ref{un_1} $ \Leftrightarrow $ \ref{un_2}) is proved.

In view of the equivalent statements  
\begin{align*}
 	\mathcal  R(A) =\mathcal  R(A^*)&\iff   \mathcal  R(A)^\perp =\mathcal  R(A^*)^\perp\\
 	&\iff  \mathcal  N(A^\dag)=\mathcal  N(A)\\
 	&\iff   \mathcal  N(A^*)=\mathcal  N(A),
 	\end{align*}
we obtain the   equivalence  (\ref{un_1} $ \iff  $ \ref{un_3} $ \iff $ \ref{un_4}).

Taking the orthogonal complements,  we obtain the implication   (\ref{un_4} $ \iff $ \ref{un_7}).  We have already observed that $\overline{C(A)}=\mathcal  N(A)^\perp$. Hence, we obtain  (\ref{un_7} $ \iff  $ \ref{un_5}).

Next we assume  (\ref{un_6}), that is,  that  $ \mathcal  H= \mathcal  R(A)\oplus \mathcal  N(A) $. Since  $ \mathcal  H=\mathcal R(A^*)  \oplus \mathcal  R(A^*)^\perp  =\mathcal R(A^*)  \oplus \mathcal  N(A)$,  we get  $\mathcal  R(A) =\mathcal  R(A^*) $.  	Now, we  assume (1), that is,    $ \mathcal  R(A) =\mathcal  R(A^*) $. Then
$$\mathcal  H=\mathcal  R(A)\oplus \mathcal  R(A)^\perp=\mathcal  R(A)\oplus \mathcal  R(A^*)^\perp=\mathcal  R(A)\oplus \mathcal  N(A)$$ 
so that   (\ref{un_6}) holds.
 \end{proof}
 
\begin{remark}
	It is proved in Theorem $\ref{ch_6_th_4}$ that 	$AA^\dag=A^\dag A$ on $ \mathcal  D(A)$ if and only if $\mathcal  N(A)=\mathcal  N(A^*)$. If we drop the assumption that $ \mathcal  R(A) $ is closed, we get that  $ AA^\dag\subset A^\dag A $ if and only if $ \mathcal  N(A)=\mathcal  N(A^*) $  and $ \mathcal  D(A^\dag)\subseteq \mathcal  D(A) $\cite{groetsch}. Indeed, if $ AA^\dag \subset A^\dag A $, then $ \mathcal  D(A^\dag)\subseteq \mathcal  D(A) $ and hence $AA^\dag =A^\dag A ~~\mbox{ on } \mathcal  D(A^\dag)$.  So $
	\overline{\mathcal  R(A)}= \overline{\mathcal  R(A^\dag)}$, and thus $\mathcal  N(A^*)=\mathcal  N(A).$	 	Conversely, if $\mathcal  N(A)=\mathcal  N(A^*)  $ and $ \mathcal  D(A^\dag)\subseteq \mathcal  D(A) $, then $\mathcal  N(A)=\mathcal  N(A^\dag)$, and hence $\overline{\mathcal  R(A^\dag)}=\overline{\mathcal  R(A)}.$ Then by Theorem \ref{thm_unpen}, we have $ AA^\dag x =A^\dag A x$ for all $ x\in \mathcal  D(A^\dag) $. Hence $ AA^\dag\subset A^\dag A  $.
		Similarly, we can prove that $A^\dag A \subset  AA^\dag$ if and only if $ \mathcal  N(A)=\mathcal  N(A^*) $  and $ \mathcal  D(A)\subseteq \mathcal  D(A^\dag) $.
\end{remark}
 To have an analogous result of Theorem \ref{bdd-EP} for closed operators, we give the following definition of $EP$ for closed operators so that Theorem $\ref{ch_6_th_4}$ will give necessary and sufficient conditions for an operator in $\mathcal  C(\mathcal  H) $ to be $EP$.
  \begin{definition}\label{defn_1} Let $A$ be a densely defined closed operator on a Hilbert space $\mathcal H$. The  operator $A$ is said to be an \textbf{$EP$ operator}  if $A$ has a closed range and  $ \mathcal  R(A)=\mathcal  R(A^*) $.  	
 \end{definition}
		
Clearly, the class of unbounded $EP$ operators include all densely defined self-adjoint operators with closed ranges.  Also every densely defined normal operator of closed range is $EP$.  Indeed, let $ A\in \mathcal  C(\mathcal  H) $  be normal with a closed range.  Then $AA^*=A^*A$ and by Theorem \ref{tarcsay}(i), we have $\mathcal D(A)=\mathcal D(A^*)$. Since $\mathcal R(A)$ is closed, $\mathcal R(A)=\mathcal R(AA^*)$ and $\mathcal R(A^*)=\mathcal R(A^*A)$. Thus $A$ is $EP$.

We now see some examples in the class of $EP$ operators.
	
 \begin{example}
Define $A$ on $\ell_2$ by 
$$A(x_1, x_2, x_3, \ldots)=(x_1, 2x_2, 3x_3, \ldots)$$ with domain 
$\mathcal D(A)= \{(x_1, x_2, x_3, \ldots)\in \ell_2 : \sum_{n=1}^\infty |n x_n|^2<\infty\}.$  
As the sequence $(1/n)_{n=1}^\infty$ belongs to $\ell_2$ but not in $\mathcal D(A)$ and $\mathcal D(A)$ contains the space $c_{00}$ of all finitely non-zero sequences, $\mathcal D(A)$ is a proper dense subsapce of $\ell_2$.  Also $\mathcal R(A)$ is closed because $A$ is surjective.  Moreover, $A$ is self-adjoint and unbounded.  Thus $A\in \mathcal C(\mathcal H)$  and it is an $EP$ operator.
 \end{example}
We have now seen an example of unbounded densely defined closed operator which is $EP$.  It is known that for $A\in \mathcal B(\mathcal H)$, $A$ is $EP$ if and only if $A^\dagger$ is $EP$, which is not true in general for $A\in \mathcal C(\mathcal H)$.  Even if $A\in \mathcal C(\mathcal H)$ has a closed range, $A^\dagger$ may not have a closed range by Theorem \ref{prop-dag} (\ref{bdd-iff-dag-closed}). In the following example, we give a densely defined closed operator satisfying the `range condition' but it is not $EP$.

 \begin{example}   Define $A$ on $\ell_2$ by 
	$$A(x_1, x_2, x_3, \ldots)=\Big(x_1, 2x_2, \frac{x_3}{3},4x_4, \frac{x_5}{5}, \ldots\Big)$$ with domain 
	$\mathcal D(A)= \{(x_1, x_2, x_3, \ldots)\in \ell_2 : (x_1, 2x_2, \frac{x_3}{3},4x_4, \frac{x_5}{5}, \ldots)\in \ell_2\}.$ 
	As the sequence $(1,\frac{1}{2},0,\frac{1}{4}, 0, \ldots)\in \ell_2\setminus\mathcal D(A)$ and $\mathcal D(A)$ contains the space $c_{00}$, $\mathcal D(A)$ is a proper dense subspace of $\ell_2$.  Since $A$ is self-adjoint, $A$ is a closed operator and 
	$\mathcal R(A)=\mathcal R(A^*)$.  Moreover, the sequence $(1/n)_{n=1}^\infty\in\ell_2\setminus\mathcal R(A)$ and $\mathcal R(A)$ is dense in $\ell^2$ as $c_{00}\subseteq \mathcal R(A)$.  Hence $\mathcal R(A)$ is a proper dense subspace of $\mathcal H$, so it cannot be closed.   Thus $A$ is not an $EP$ operator.	
\end{example}

 \begin{example} \cite{Huang2012}
 	Let   $ \varphi :[0,1]\rightarrow \mathbb{C} $ by $$\varphi  (t)=\left\{\begin{array}{ccl}
 	1 & \text{  if } & t=0\\
 	\frac{1}{\sqrt{t}} & \text{  if } & 0<t\leq 1.	\end{array}\right.$$
 	Define $$ Af=\varphi f,\quad  f\in \mathcal D(A), $$ 
 	where 
 	$ \mathcal D(A)=\left\{f\in L^2[0,1]: \varphi f\in L^2[0,1]\right\}.$  Then $ A $ is  a densely defined closed operator. As $ |\varphi (t)|\geq 1 $ for all $ t\in [0,1] $, we have $ \mathcal R(A)=L^2[0,1] $ and $ A $ has bounded inverse $ A^{-1} :L^2[0,1]\rightarrow L^2[0,1] $ defined by $ A^{-1}g  = \psi g $ for all $ g\in L^2[0,1] $
 	where 
 	$$\psi (t)=\left\{\begin{array}{ccl}
 	1 & \text{  if } & t=0\\
 	\sqrt{t}& \text{  if } & 0<t\leq 1.	\end{array}\right.$$
 	Hence $ A $ is an $ EP $ operator on $L^2[0,1]  $.
 \end{example} 
 
 \begin{example} 
 	Let $\mathcal H=L^2[0,1]$. Let
 	$$\mathcal{AC}[0,1]=\Big\{f\in \mathcal H: f:[0,1]\rightarrow \mathbb C \ \  \text{is absolutely continuous and}\; f'\in\mathcal H \Big\}.$$ Let
 	$\mathcal D(A)=\Big\{f\in \mathcal {AC}[0,1]: f(0)=f(1) \Big\}$. 	Define  $A:\mathcal D(A)\rightarrow \mathcal H$ by $$Af=if'\quad \text{for all}\; f\in
 	\mathcal D(A).$$ 
 	It can be shown that $A$ is a densely defined self-adjoint operator and $$\mathcal R(A)=\Big\{u\in \mathcal H:
 	\int\limits_0^1 u(t){d}t=0 \Big\}={span\{1}\}^\bot.$$ Hence $A$ is an $EP$ operator.
 \end{example} 
 
Now, we prove some results related to unbounded $EP$ operators.

		\begin{proposition}
			Let $ A\in \mathcal  C(\mathcal  H) $ with a closed range. If $ A $ is $ EP $, then the following statements are true.
			\begin{enumerate}
				\item $A^*$ is $EP$.\label{prop1.2}
				\item $AA^*$ is $EP$.\label{prop1.3}
				\item $A^*A$ is $EP$. \label{prop1.4}
				\item $|A|$ is  $EP$. \label{1.5}
				
			\end{enumerate}
		\end{proposition}
		\begin{proof}
			The  proof of (\ref{prop1.2}) is obvious from the definition of $ EP $ operator and  Theorem \ref{Ben_thm_1}. Since $ AA^* $ and $ A^*A $ are self-adjoint and by Theorem \ref{Ben_thm_1}, $ AA^* $ and $ A^*A $ are $ EP $. 	Since $ \mathcal  R(A) $ is closed, by  Theorem \ref{Ben_thm_1}, we have $ \mathcal  R(A^*) = \mathcal  R(A^*A)$. Also, we have $\overline{\mathcal  R(|A|)}=\overline{\mathcal  R(A^*)}$. Now $ \mathcal  R(A^*A)\subseteq \mathcal  R(|A|) \subseteq \overline{\mathcal  R(|A|)}=\overline{\mathcal  R(A^*)} =\mathcal  R(A^*) =\mathcal  R(A^*A) $. Hence $ \mathcal  R(A^*)=\mathcal  R(|A|) $. Therefore $ \mathcal  R(|A|) $ is closed. Since $ |A| $ is self-adjoint and $ \mathcal  R(|A|) $ is closed, $ |A| $ is $ EP $.	
		\end{proof}

  \begin{theorem}
 	Let $ A\in \mathcal  C(\mathcal  H) $ with a closed range. Then the following are equivalent:
 	\begin{enumerate}
 		\item $A$  is $EP$ ;
 		\item $AA^\dag=A^\dag A$ on $ \mathcal  D(A) $ ;
 		\item $\mathcal  N(A)=\mathcal  N(A^\dag)$ ;
 		\item $\mathcal  N(A)=\mathcal  N(A^*)$ ;
 		\item $\mathcal  N(A)^\perp =\mathcal  R(A)$ ; 
 		\item  $ \overline{C(A)}=\mathcal  R(A) $ ;
 		\item $ \mathcal  H= \mathcal  R(A)\oplus \mathcal  N(A)$.
 	\end{enumerate}
 \end{theorem}
 \begin{proof} Follows from Theorem \ref{ch_6_th_4}.
 \end{proof}
		\begin{remark}
	If $ A\in \mathcal  C(\mathcal  H) $ is an unbounded EP operator, then $\mathcal  R(A^\dag)$ is a proper dense subspace of $\mathcal  R(A)$ because $\overline{\mathcal R(A^\dag)}=\overline{ \mathcal  R(A^*)}=\mathcal R(A^*)=\mathcal R(A)$ and by Theorem \ref{prop-dag} (\ref{bdd-iff-dag-closed}).
	
\end{remark}

 \begin{theorem}\label{extra_1}
 	Let $ A\in \mathcal  C(\mathcal  H) $ be $EP$.  Then there exists a densely defined bijective linear operator  $C:\mathcal D(A^*)\to \mathcal D(A)$  such that  $A^*\subset AC$. 
 \end{theorem}

 \begin{proof}
Assume that $A$ is $EP$. Let $ y\in \mathcal  D(A^*) $. Then $y=y_1+y_2, y_1\in C(A^*), y_2\in \mathcal  N(A^*)$.	  Since $\mathcal R(A^*)=\mathcal R(A)$, $A^*y=Ax$, for some $ x\in \mathcal  D(A)$ with $ x=x_1+x_2, x_1\in C(A), x_2\in \mathcal  N(A)$. Define $ Cy=C(y_1+y_2)=	y_1+x_2 $.  It is easily seen that $C$ is well-defined, linear and $A^*y=ACy$, for all $y\in \mathcal D(A^*)$.   
	
	We shall now prove $C$ is injective. Suppose that $C(y_1+y_2)=0$, so $y_1+x_2=0$, hence $y_1=x_2=0$ because $\mathcal N(A^*)=\mathcal N(A)$.  Also, $A^*y_2=0$.  Therefore $y_2=0$, thus $y_1+y_2=0.$  To prove surjective, let us take $x=x_1+x_2\in \mathcal D(A)$.  Then there exists $y\in y_1+y_2\in \mathcal D(A^*)$  such that $A^*y=Ax$.  Let $z=x_1+y_2$.  Then $Cz=x$, hence $C$ is bijective.
	
 \end{proof}

 \section{Closed Hypo-$ EP $ Operators on Hilbert Spaces}
 Weakening the range-space condition by inclusion relation between range spaces, the notion of hypo-$EP$ was defined by Itoh for bounded operators in \cite{Itoh2005}. Algebraic and analytic characterizations, sum, product, factorization of bounded hypo-$EP$ operators are discussed in \cite{Sam2017, Sam2018}.  In the collection of bounded operators, the class of all hypo-$EP$ operators contains the class of all $EP$ operators and hyponormal operators with closed ranges. Hence it contains all normal and invertible operators with closed ranges.  In the case of finite dimensional, $EP$ and hypo-$EP$ are same.   In this section, we define and discuss hypo-$EP$ operator for densely defined closed operators. 
 \begin{definition}  Let $A$ be a densely defined closed operator on a Hilbert space $\mathcal H$. The  operator $A$ is said to be a \textbf{hypo-$ EP $ operator} if $A$ has a closed range and $ \mathcal  R(A)\subseteq \mathcal  R(A^*) $.
 \end{definition}
 	In the above definition, the inclusion relation for range spaces ``$ \mathcal  R(A)\subseteq \mathcal  R(A^*) $'' can equivalently be replaced by inclusion relation for null spaces $ ``\mathcal  N(A)\subseteq \mathcal  N(A^*) $''.

 \begin{example}
 	Define $A$ on $\ell_2$ by  $$A(x_1,x_2,x_3,\ldots)=(0,x_1,2x_2,3x_3,\ldots) $$ with 	$\mathcal  D(A)=\left\{(x_1,x_2,x_3,\ldots)\in \mathcal  H: 
 	\sum_{n=1}^\infty |nx_n|^2<\infty \right\}.$
As the sequence $(1/n)_{n=1}^\infty$ belongs to $\ell_2$ but not in $\mathcal D(A)$ and $\mathcal D(A)$ contains the space $c_{00}$, $\mathcal D(A)$ is a proper dense subsapce of $\ell_2$.  Also $A$ is closed and $\mathcal R(A)=\ell_2\setminus span\{e_2\}$. The adjoint of $A$ is defined by $A^*(x_1, x_2, x_3, \ldots)=(x_2, 2x_3,3x_4, \ldots)$ with $\mathcal D(A^*)=\left\{(x_1,x_2,x_3,\ldots)\in \mathcal  H: \sum_{n=2}^\infty |(n-1)x_n|^2<\infty \right\}.$  Moreover, $\mathcal N(A^*)=span\{e_1\}$ and $\mathcal N(A)=\{0\}$, so $ \mathcal  R(A)\subseteq \mathcal  R(A^*) $.  Thus $A$ is hypo-$EP$ but not $EP$.
	
 \end{example}
 \begin{theorem}\label{hypo-first-char}
 	Let $ A\in \mathcal  C(\mathcal  H) $ with a closed range. Then $ A $ is hypo-$ EP $ if and only if $ A^\dag A^2A^\dag=AA^\dag $.
 	
 \end{theorem}
 \begin{proof}
 	Suppose $ \mathcal  R(A) \subseteq \mathcal  R(A^*)$ and $ \mathcal  R(A) $ is closed. Then $ AA^\dag x\in \mathcal  R(A) $  for each $ x\in \mathcal  H $ and hence $ AA^\dag x\in \overline{\mathcal  R(A^\dag)}=\mathcal  R(A^*) $. As $ A^\dag A $ is a projection onto $ \overline{\mathcal  R(A^\dag)} $, we have $ A^\dag A(AA^\dag x)=AA^\dag x $. Hence $ A^\dag A^2A^\dag=AA^\dag $. Conversely, suppose $ A^\dag A^2A^\dag=AA^\dag $. Then  $\mathcal  R(A)=\overline{\mathcal  R(A)}=\mathcal  R(AA^\dag)=\mathcal  R(A^
 	\dag A^2A^\dag)\subseteq \mathcal  R(A^\dag)=C(A)\subseteq \mathcal  N(A)^\perp=\overline{\mathcal  R(A^*)}=\mathcal  R(A^*)  $. Hence $ A $ is hypo-$EP$. \end{proof} \begin{theorem}\label{thunhy_1}
 	Let $ A\in \mathcal  C(\mathcal  H) $. Then each of the following statements implies the next statement:
 	\begin{enumerate}
 		\item $ A $ is hypo-$ EP $ ;\label{unhy_1}
 		\item $ A(A^\dag)^2A=AA^\dag $ on $ \mathcal  D(A) $ ; \label{unhy_2}
 		\item $ AA^\dag \leq A^\dag A $ on $ \mathcal  D(A) $ ; \label{unhy_3}
 		\item $ \|AA^\dag x\|\leq \|A^\dag Ax\| $ for all $ x\in \mathcal  D(A) $.\label{unhy_4}
 	\end{enumerate}
 \end{theorem}
 \begin{proof}
 	Assume that  $ A $ is hypo-$ EP $  and  $ x\in \mathcal  D(A) $. Then 
 	\begin{align*}
 	\left<AA^\dag A^\dag Ax,x\right>&=\left<(AA^\dag)^*A^\dag Ax,x\right>\\
 	&=\left<A^\dag Ax, AA^\dag x\right>\\
 	&=\left<(A^\dag A)^*x, AA^\dag x\right>\\
 	&= \left<x, A^\dag A^2A^\dag x\right>\\
 	&= \left<x, AA^\dag x\right>, \text{by Theorem} \ \ref{hypo-first-char}\\
 	&= \left<AA^\dag x, x\right>.
 	\end{align*}
 	Hence $ A(A^\dag)^2A=AA^\dag $ on $ \mathcal  D(A) $. 
 	
 	\noindent Assume  that $ A(A^\dag)^2A=AA^\dag $ on $ \mathcal  D(A) $. Let  $ x\in \mathcal  D(A) $. Then
 	\begin{align*}
 	\left<AA^\dag x,x\right>&=\left<AA^\dag  AA^\dag x, x\right>\\
 	&=\left<(AA^\dag )^* AA^\dag x, x\right>\\
 	&=\left< AA^\dag x,(AA^\dag ) x\right>\\
 	&=\|AA^\dag x\|^2\\
 	&=\|A(A^\dag)^2Ax\|^2\\
 	&\leq \|AA^\dag\|^2\|A^\dag Ax\|^2\\
 	&=\|A^\dag Ax\|^2\\
 	&=\left<A^\dag Ax,A^\dag Ax\right>\\
 	&=\left<A^\dag Ax,x\right>.
 	\end{align*}
 	Hence $ AA^\dag \leq A^\dag A $ on $ \mathcal  D(A) $. 
 	
 	\noindent Assume that  $ AA^\dag \leq A^\dag A $ on $ \mathcal  D(A) $. Let $ x\in \mathcal  D(A) $. Then 
 	\begin{align*}
 	&\left<AA^\dag x,x \right>\leq \left<A^\dag Ax,x\right>\\
 	\Rightarrow & \left<AA^\dag AA^\dag x,x \right>\leq \left<A^\dag AA^\dag Ax,x\right>\\
 	\Rightarrow & \left< AA^\dag x, AA^\dag x \right>\leq \left<A^\dag Ax,A^\dag Ax\right>\\
 	\Rightarrow & \|AA^\dag x\|^2\leq \|A^\dag Ax\|^2.
 	\end{align*}
 	Thus	$ \|AA^\dag x\leq \|A^\dag Ax\| $ for all $ x\in \mathcal  D(A)$.	
 \end{proof}

It is observed that $A$ is hypo-$EP$ if and only if $AA^\dagger=A^\dagger A$, for $A\in \mathcal B(\mathcal H)$.  This is not true for $A\in \mathcal C(\mathcal H)$.  However, if $ \mathcal  R(A) \subseteq \mathcal  D(A)$, all the necessary conditions for hypo-$ EP $ in Theorem \ref{thunhy_1} become sufficient conditions as well. Note that the inclusion relation ``$ \mathcal  R(A) \subseteq \mathcal  D(A)$'' is irredutant in the case of bounded operators.
 
 \begin{theorem}\label{th_u1}
 	Let $A\in\mathcal  C(\mathcal  H)$ be hypo-$EP$.  Then for each $x\in \mathcal D(A)$, there exists $k>0$ such that $ 	|\langle Ax, y\rangle|\leq k \|Ay\|,\mbox{ for all }  y\in \mathcal  D(A). $
 	
 \end{theorem}
 
 \begin{proof}
 	Suppose $A$ is hypo-$EP$.  If $x\in \mathcal  N(A)$, then the result is trivial.  Let $x\in \mathcal  D(A)$ such that $Ax\neq 0$. Then $Ax\in \mathcal {R}(A)\subseteq\mathcal {R}(A^*)$. Therefore there exists   $z\in \mathcal  D(A^*)$ such that $A^*z=Ax$. Then for all $y\in \mathcal  D(A)$,
 	\begin{eqnarray*}
 		|\langle Ax,y\rangle|=|\langle A^*z,y\rangle|=|\langle z,Ay\rangle| \leq\|z\| \|Ay\|.
 	\end{eqnarray*}
 	Taking $k=\|z\|$, we get $$|\langle Ax, y\rangle|\leq k \|Ay\|,$$ for all $y\in \mathcal  D(A)$.	
 \end{proof}
 
 The converse of Theorem \ref{th_u1} has been proved for bounded hypo-$ EP $ operators on Hilbert spaces in \cite{Sam2017} where Douglas' theorem for bounded operators was used. Unlike the bounded operators, Douglas' theorem for densely defined closed operators does not guarantee the equivalance of the  notions of majorization, range inclusion and factorization.

\section{A perturbation result}

Stability of $EP$ operator under perturbation by bounded operators is investigated in the following result.

 \begin{theorem}
 	Let $ A\in \mathcal  C(\mathcal  H) $  be an   $EP$ operator. Let $B\in \mathcal  B(\mathcal  H)$ be such that $\|B\| \|A^\dagger\|<1,  BA^\dagger A=B|_{\mathcal  D(A)}$ and $AA^\dagger B=B$. Then $A+B$ is $EP$. 
 \end{theorem}
 \begin{proof}
 	Let $ x\in \mathcal  N(A) $. Then $ Bx=BA^\dag Ax=0 $. Hence $ (A+B)x=Ax+Bx=0 $ and $ x\in \mathcal  N(A+B) $. Therefore $ \mathcal  N(A)\subseteq \mathcal  N(A+B)$. It is given that $\|B\| \|A^\dagger\|<1$, hence $ \|BA^\dag\|<1 $ which implies $ I+BA^\dag $ is invertible. Let $ x\in \mathcal  N(A+B) $. Then
 	\begin{eqnarray*}
 		(A+B)x&=&0\\
 		Ax+BA^\dag Ax&=&0\\
 		(I+BA^\dag)Ax&=&0\\
 		Ax&=&0.
 	\end{eqnarray*}
 	Therefore $ x\in \mathcal  N(A) $ and $ \mathcal  N(A+B)\subseteq \mathcal  N(A)$. Hence $ \mathcal  N(A+B)=\mathcal  N(A)$.
 	
 	Since $ \mathcal  R(A) $ is closed,  $ \|Ax\|\geq \gamma (A) \|x\| $, for all $ x\in C(A) $. Since $ \mathcal  N(A+B)=\mathcal  N(A) $, we have $ C(A+B)=C(A) $. For any $ x\in C(A+B) $, $ \|(A+B)x\|=\|Ax+Bx\|\geq \Big| \|Ax\|-\|Bx\| \Big| \geq \gamma(A) \|x\| - \|B\| \|x\|= (\gamma (A)-\|B\|)\|x\|$. Since $ \gamma(A)=\frac{1}{\|A^\dag\|}$ and $\|B\| \|A^\dagger\|<1$, we have $ \gamma(A)-\|B\|>0 $. Hence $ \mathcal  R(A+B) $ is closed.
 	
 	Let $ y\in \mathcal  R(A+B) $. Then there exists $ x\in \mathcal  D(A) $ such that $ y=(A+B)x=Ax+Bx=Ax+AA^\dag B x=A(I+A^\dag B) x$. Hence $ y\in \mathcal  R(A(I+A^\dag B))\subseteq \mathcal  R(A) $  and $ R(A+B)\subseteq \mathcal  R(A) $.  Let $ y\in \mathcal  R(A) $. Then there exists $ x\in \mathcal  D(A) $ such that $ Ax=y $.  Since $\|B\| \|A^\dagger\|<1$, $ \|A^\dag B\|<1 $ implies $ (I+A^\dag B)^{-1}\in \mathcal  B(\mathcal  H) $. Since $ (I+A^\dag B) $ is surjective, for $ x\in \mathcal  H $, there exists $ x^\prime \in \mathcal  D(A) $ such that $ (I+A^\dag B)x^\prime  =x $. Hence $ y=Ax=A(I+A^\dag B)x^\prime= Ax^\prime +AA^\dag Bx^\prime = Ax^\prime +Bx^\prime=(A+B)x^\prime $. Therefore $ y\in \mathcal  R(A+B) $ and hence $ \mathcal  R(A)=\mathcal  R(A+B) $.
 	
 	Since $A$ is $EP$, $\mathcal  R(A)=\mathcal  R(A^*)$.  Since $\mathcal  N(A+B)=\mathcal  N(A)$, $\mathcal  R(A^*+B^*)=\mathcal  R(A^*).$   Thus $\mathcal  R(A+B)=\mathcal  R(A)=\mathcal  R(A^*)=\mathcal  R(A^*+B^*)$, so $A+B$ is $EP$.  	
 \end{proof}
 
  \section{Acknowledgement}
  
  The author is thankful to Prof. M. Thamban Nair for having some discussion while preparing the manuscript and also for providing the proof for Theorem \ref{extra_1}.


\begin{thebibliography}{10}
 	
 	\bibitem{Ben-Israel2003}
 	Adi Ben-Israel and Thomas N.~E. Greville.
 	\newblock {\em Generalized inverses}, volume~15 of {\em CMS Books in
 		Mathematics/Ouvrages de Math\'{e}matiques de la SMC}.
 	\newblock Springer-Verlag, New York, second edition, 2003.
 	\newblock Theory and applications.
 	
 	\bibitem{MR0226433}
 	S.~J. Bernau.
 	\newblock The square root of a positive self-adjoint operator.
 	\newblock {\em J. Austral. Math. Soc.}, 8:17--36, 1968.
 	
 	\bibitem{Brock1990}
 	K.~G. Brock.
 	\newblock A note on commutativity of a linear operator and its
 	{M}oore-{P}enrose inverse.
 	\newblock {\em Numer. Funct. Anal. Optim.}, 11(7-8):673--678, 1990.
 	
 	\bibitem{Campbell1975}
 	Stephen~L. Campbell and Carl~D. Meyer.
 	\newblock {{$EP$} operators and generalized inverses}.
 	\newblock {\em Canad. Math. Bull}, 18(3):327--333, 1975.
 	
 	\bibitem{Koliha2007}
 	Dragan~S. Djordjevi{\'c} and J.~J. Koliha.
 	\newblock {Characterizing {H}ermitian, normal and {EP} operators}.
 	\newblock {\em Filomat}, 21(1):39--54, 2007.
 	
 	\bibitem{Douglas1966}
 	R.~G. Douglas.
 	\newblock On majorization, factorization, and range inclusion of operators on
 	{H}ilbert space.
 	\newblock {\em Proc. Amer. Math. Soc.}, 17:413--415, 1966.
 	
 	\bibitem{groetsch}
 	C.~W. Groetsch.
 	\newblock Inclusions and identities for the {M}oore-{P}enrose inverse of a
 	closed linear operator.
 	\newblock {\em Math. Nachr.}, 171:157--164, 1995.
 	
 	\bibitem{Hartwig1997}
 	Robert~E. Hartwig and Irving~J. Katz.
 	\newblock {On products of {EP} matrices}.
 	\newblock {\em Linear Algebra Appl.}, 252:339--345, 1997.
 	
 	\bibitem{Huang2012}
 	Qianglian Huang, Lanping Zhu, and Jiena Yu.
 	\newblock Some new perturbation results for generalized inverses of closed
 	linear operators in {B}anach spaces.
 	\newblock {\em Banach J. Math. Anal.}, 6(2):58--68, 2012.
 	
 	\bibitem{Itoh2005}
 	Masuo Itoh.
 	\newblock On some {EP} operators.
 	\newblock {\em Nihonkai Math. J.}, 16(1):49--56, 2005.
 	
 	\bibitem{Pearl1966}
 	Irving~Jack Katz and Martin~H. Pearl.
 	\newblock On {$EPr$} and normal {$EPr$} matrices.
 	\newblock {\em J. Res. Nat. Bur. Standards Sect. B}, 70B:47--77, 1966.
 	
 	\bibitem{nair-linop}
 	M.~Thamban Nair.
 	\newblock {\em Linear Operator Equations: Approximations and Regularization}.
 	\newblock World Scientific, first edition, 2009.
 	
 	\bibitem{nair-fa}
 	M.~Thamban Nair.
 	\newblock {\em Functional Analysis: A First Course}.
 	\newblock Prentice-Hall of India, second edition, 2021.
 	
 	\bibitem{Nashed1976}
 	M.~Zuhair Nashed, editor.
 	\newblock {\em Generalized inverses and applications}.
 	\newblock Academic Press [Harcourt Brace Jovanovich, Publishers], New
 	York-London, 1976.
 	\newblock University of Wisconsin, Mathematics Research Center, Publication No.
 	32.
 	
 	\bibitem{Pearl1966(2)}
 	M.~H. Pearl.
 	\newblock On generalized inverses of matrices.
 	\newblock {\em Proc. Cambridge Philos. Soc.}, 62:673--677, 1966.
 	
 	\bibitem{Rudin1991}
 	Walter Rudin.
 	\newblock {\em Functional analysis}.
 	\newblock International Series in Pure and Applied Mathematics. McGraw-Hill,
 	Inc., New York, second edition, 1991.
 	
 	\bibitem{Sam2018}
 	P.~Sam~Johnson and A.~Vinoth.
 	\newblock Product and factorization of hypo-{EP} operators.
 	\newblock {\em Spec. Matrices}, 6:376--382, 2018.
 	
 	\bibitem{sch-book}
 	Konrad Schm\"{u}dgen.
 	\newblock {\em Unbounded self-adjoint operators on {H}ilbert space}, volume 265
 	of {\em Graduate Texts in Mathematics}.
 	\newblock Springer, Dordrecht, 2012.
 	
 	\bibitem{Hans1950}
 	Hans Schwerdtfeger.
 	\newblock {\em Introduction to {L}inear {A}lgebra and the {T}heory of
 		{M}atrices}.
 	\newblock P. Noordhoff, Groningen, 1950.
 	
 	\bibitem{tarcsay-2012}
 	Zs. Tarcsay.
 	\newblock Operator extensions with closed range.
 	\newblock {\em Acta Math. Hungar.}, 135(4):325--341, 2012.
 	
 	\bibitem{Sam2017}
 	A.~Vinoth and P.~Sam~Johnson.
 	\newblock On sum and restriction of hypo-{$EP$} operators.
 	\newblock {\em Funct. Anal. Approx. Comput.}, 9(1):37--41, 2017.
 	
 	\bibitem{yosida}
 	K\^{o}saku Yosida.
 	\newblock {\em Functional analysis}, volume 123 of {\em Grundlehren der
 		Mathematischen Wissenschaften [Fundamental Principles of Mathematical
 		Sciences]}.
 	\newblock Springer-Verlag, Berlin-New York, sixth edition, 1980.
 	
 \end{thebibliography}
\end{document}